# General Design Bayesian Generalized Linear Mixed Models

**Y. Zhao, J. Staudenmayer, B. A. Coull and M. P. Wand**

*Abstract.* Linear mixed models are able to handle an extraordinary range of complications in regression-type analyses. Their most common use is to account for within-subject correlation in longitudinal data analysis. They are also the standard vehicle for smoothing spatial count data. However, when treated in full generality, mixed models can also handle spline-type smoothing and closely approximate kriging. This allows for nonparametric regression models (e.g., additive models and varying coefficient models) to be handled within the mixed model framework. The key is to allow the random effects design matrix to have general structure; hence our label *general design*. For continuous response data, particularly when Gaussianity of the response is reasonably assumed, computation is now quite mature and supported by the R, SAS and S-PLUS packages. Such is not the case for binary and count responses, where generalized linear mixed models (GLMMs) are required, but are hindered by the presence of intractable multivariate integrals. Software known to us supports special cases of the GLMM (e.g., PROC NLMIXED in SAS or glmmML in R) or relies on the sometimes crude Laplace-type approximation of integrals (e.g., the SAS macro glimmix or glmmPQL in R). This paper describes the fitting of general design generalized linear mixed models. A Bayesian approach is taken and Markov chain Monte Carlo (MCMC) is used for estimation and inference. In this generalized setting, MCMC requires sampling from nonstandard distributions. In this article, we demonstrate that the MCMC package WinBUGS facilitates sound fitting of general design Bayesian generalized linear mixed models in practice.

*Y. Zhao is Mathematical Statistician, Division of Biostatistics, Center for Devices and Radiological Health, U.S. Food and Drug Administration, 1350 Piccard Drive, Rockville, Maryland 20850, USA (e-mail: yihua.zhao@fda.hhs.gov). J. Staudenmayer is Assistant Professor, Department of Mathematics and Statistics, University of Massachusetts, Amherst, Massachusetts 01003, USA (e-mail: jstauden@math.umass.edu). B. A. Coull is Associate Professor, Department of Biostatistics, School of Public Health, Harvard University, 665 Huntington Avenue, Boston, Massachusetts 02115, USA (e-mail:*

*bcoull@hsph.harvard.edu). M. P. Wand is Professor, Department of Statistics, School of Mathematics, University of New South Wales, Sydney 2052, Australia (e-mail: wand@maths.unsw.edu.au).*

## 1. INTRODUCTION

The generalized linear mixed model (GLMM) is one of the most useful structures in modern statistics, allowing many complications to be handled within the familiar linear model framework. The fitting of such models has been the subject of a great deal of research over the past decade. Early contributions to fitting various forms of the GLMM include Stiratelli, Laird and Ware (1984), Anderson and Aitkin (1985), Gilmour, Anderson and Rae (1985), Schall (1991), Breslow and Clayton (1993) and Wolfinger and O'Connell (1993). A summary is provided by McCulloch and Searle (2001, Chapter 10).

Most of the literature on fitting GLMMs is geared toward grouped data. Examples include repeated binary responses on a set of subjects and standardized mortality ratios in geographical subregions. However, GLMMs are much richer than the subclass needed for these situations. The key to full generality is the use of *general design matrices*, for both the fixed and random components. Once again, we refer to McCulloch and Searle (2001, Chapter 8) for an overview of general design GLMMs. An excellent synopsis of general design linear mixed models is provided by Robinson (1991) and the ensuing discussion. One of the biggest payoffs from the general design framework is the incorporation of nonparametric regression, or smoothing, through penalized regression splines (e.g., Wahba, 1990; Speed, 1991; Verbyla, 1994; Brumback, Ruppert and Wand, 1999). Higher dimensional extensions essentially correspond to generalized kriging (Diggle, Tawn and Moyeed, 1998). This allows for smoothing-type models such as generalized additive models to be fitted as a GLMM and combined with the more traditional grouped data uses. This is the main thrust of the recent book by Ruppert, Wand and Carroll (2003), a summary of which is provided by Wand (2003). General designs also permit the handling of crossed random effects (e.g., Shun, 1997) and multilevel models (e.g., Goldstein, 1995; Kreft and de Leeuw, 1998).

The simplest method for fitting general design GLMMs involves Laplace approximation of integrals (Breslow and Clayton, 1993; Wolfinger and O'Connell, 1993) and is commonly referred to as *penalized quasi-likelihood* (PQL). However, the approximation can be quite inaccurate in certain circumstances. Breslow and Lin (1995) and Lin and Breslow (1996) showed that PQL leads to estimators that are asymptotically biased. For situations such as paired binary data the PQL approximation is particularly poor. In their summary of PQL, McCulloch and Searle (2001, Chapter 10, pages 283–284) concluded by stating that they "cannot recommend the use of simple PQL methods in practice." In this article we take a Bayesian approach and explore the Markov chain Monte Carlo (MCMC) fitting of general design GLMMs. One advantage of a Bayesian approach over its frequentist counterpart includes the fact that uncertainty in variance components is more easily taken into account (e.g., Handcock and Stein, 1993; Diggle, Tawn and Moyeed, 1998). As summarized in Section 9.6 of McCulloch and Searle (2001), the frequentist approach to this problem is thwarted by largely intractable distribution theory. Under a Bayesian approach, posterior distributions of parameters of interest take this variability into account. The hierarchical structure of the Bayesian GLMMs lends itself to Gibbs sampling schemes, albeit with some nonconjugate full conditionals, to sample from these posteriors. In addition, it is computationally simpler to obtain variance estimates of the predictions of the random effects. Booth and Hobert (1998) showed that, in a frequentist framework, second-order estimation of the conditional standard error of prediction for the random effects requires bootstrapping the maximum likelihood estimates of the fixed effects and variance components. For complicated random effects structures, computation of a single maximum likelihood fit can be expensive, making the bootstrap computationally prohibitive. In the Bayesian framework, interest focuses on the posterior variance of the random effects given the data, which is a by-product of the MCMC output. Note, however, that the Bayesian approach involves specification of prior distributions of all model parameters. This requires some care, especially when sample sizes are small.

There have been a few other contributions to Bayesian formulations of GLMMs in the literature. Those known to us are Zeger and Karim (1991), Clayton (1996), Diggle, Tawn and Moyeed (1998) and Fahrmeir and Lang (2001). However, each of these articles is geared toward special cases of GLMMs. The GLMMs described in this article are much more general and allow for random effects models for longitudinal data, crossed random effects, smoothing of spatial count data, generalized additive models, generalized geostatistical models, additive models with interactions, varying coefficient models and various combinations of these (Wand, 2003).

Section 2 lays out notation for general design GLMMs and gives several important examples. MCMC implementation is described in Section 3, with a focus on the `WinBUGS` package. Section 4 provides three illustratory data analyses. We close with some discussion in Section 5.

## 2. MODEL FORMULATION

The GLMMs for canonical one-parameter exponential families (e.g., Poisson, logistic) and Gaussian random effects take the general form

$$
\begin{aligned}
[\mathbf{y}|\boldsymbol{\beta},\mathbf{u}] = \exp\{ & \mathbf{y}^\top(\mathbf{X}\boldsymbol{\beta}+\mathbf{Z}\mathbf{u}) \\
& - \mathbf{1}^\top b(\mathbf{X}\boldsymbol{\beta}+\mathbf{Z}\mathbf{u}) + \mathbf{1}^\top c(\mathbf{y})\},
\end{aligned}
\tag{1}
$$

$$
[\mathbf{u}|\mathbf{G}] \sim N(\mathbf{0},\mathbf{G}),
\tag{2}
$$

where here and throughout the distribution of a random vector $\mathbf{x}$ is denoted by $[\mathbf{x}]$ and the conditional distribution of $\mathbf{y}$ given $\mathbf{x}$ is denoted by $[\mathbf{y}|\mathbf{x}]$.

In the Poisson case $b(x)=e^x$, while in the logistic case $b(x)=\log(1+e^x)$. A few other models (e.g., gamma, inverse Gaussian) also fit into this structure (McCullagh and Nelder, 1989). A number of extensions and modifications are possible. One is to allow for overdispersion, especially in the Poisson case. In this paper we will restrict attention to the canonical one-parameter exponential family structure. In most situations, the main parameters of interest are contained in $\boldsymbol{\beta}$ and $\mathbf{G}$, and prior distributions for them need to be specified; see Sections 2.1 and 2.2.

It is important to separate out random effects structure for handling grouping. One reason is that this allows for the possibility of hierarchical centering in the MCMC implementations (Section 2.3). It also recognizes the different covariance structures used in longitudinal data modeling, smoothing and spatial statistics. Such considerations suggest the breakdown

$$
\begin{aligned}
\mathbf{X}\boldsymbol{\beta}+\mathbf{Z}\mathbf{u} = {} & \mathbf{X}^R\boldsymbol{\beta}^R + \mathbf{Z}^R\mathbf{u}^R \\
& + \mathbf{X}^G\boldsymbol{\beta}^G + \mathbf{Z}^G\mathbf{u}^G + \mathbf{Z}^C\mathbf{u}^C,
\end{aligned}
\tag{3}
$$

where

$$
\mathbf{X}^R \equiv \begin{bmatrix} \mathbf{X}_1^R \\ \vdots \\ \mathbf{X}_m^R \end{bmatrix}, \quad \mathbf{Z}^R \equiv \operatorname*{blockdiag}_{1\le i\le m}(\mathbf{X}_i^R)
$$

and

$$
\operatorname{Cov}(\mathbf{u}^R) \equiv \operatorname*{blockdiag}_{1\le i\le m}(\boldsymbol{\Sigma}^R) \equiv \mathbf{I}_m \otimes \boldsymbol{\Sigma}^R
$$

correspond to random intercepts and slopes, as typically used for repeated measures data on $m$ groups with sample sizes $n_1,\ldots,n_m$. Here $\mathbf{X}_i^R$ is an $n_i \times q^R$ matrix for the random design corresponding to the $i$th group, $\boldsymbol{\Sigma}^R$ is an unstructured $q^R \times q^R$ covariance matrix and $\otimes$ denotes Kronecker product.

Next, $\mathbf{X}^G$ and $\mathbf{Z}^G$ are general design matrices, usually of different form than those arising in random effects models. In many of our examples, $\mathbf{X}^G$ contains indicator variables or polynomial basis functions of a continuous predictor, while $\mathbf{Z}^G$ contains spline basis functions (e.g., Brumback, Ruppert and Wand, 1999). The $\mathbf{Z}^G\mathbf{u}^G$ term may be further decomposed as

$$
\mathbf{Z}^G\mathbf{u}^G = \sum_{\ell=1}^{L} \mathbf{Z}_\ell^G \mathbf{u}_\ell^G
$$

with each $\mathbf{Z}_\ell^G$, $1 \le \ell \le L$, usually corresponding to a smooth term in an additive model. Also, in keeping with spline penalization, we only consider

$$
\operatorname{Cov}(\mathbf{u}^G) = \operatorname*{blockdiag}_{1\le \ell \le L}(\sigma_{u\ell}^2 \mathbf{I}).
$$

Note that the decomposition (3) is not unique for a particular model. For instance, in the crossed random effects model given in the following Example 3, we present two methods of decomposition.

The $\mathbf{Z}^C\mathbf{u}^C$ component represents random effects with spatial correlation structure. This can be done in a number of ways (e.g., Wakefield, Best and Waller, 2001); we will just describe one of the more common approaches here. Suppose disease incidence data are available over $N$ contiguous regions. The random effect $\mathbf{u}^C$ vector is of dimension $N$ with entries $U_1^C$, $\ldots, U_N^C$. The conditional distribution of $U_i^C$ given $U_j^C$, $j \ne i$, is a univariate normal distribution with mean equal to the average $U_j^C$ values of $U_i^C$'s neighboring regions and variance equal to $\sigma_c^2$ divided by the number of neighboring regions. This is known as the intrinsic Gaussian autoregression distribution (Besag, York and Mollié, 1991). This leads to $\mathbf{u}^C$ having an improper density proportional to

$$
\exp\left\{ -\sum_{i \sim j} \tfrac{1}{2}\sigma_c^{-2}(U_i^C - U_j^C)^2 \right\},
\tag{4}
$$

where $i \sim j$ denotes spatially adjacent regions.

The versatility of (3) can be appreciated by considering the following set of examples. Note that we use truncated linear basis functions for smoothing components to keep the formulations simple (e.g.,

Brumback, Ruppert and Wand, [1999]). In practice these may be replaced by B-splines (Durbán and Currie, [2003]) or radial basis functions (French, Kammann and Wand, [2001]). Knots are denoted by $\kappa_k$ with possible superscripting. Ruppert ([2002]) discussed choice of knots of univariate smoothings, whereas Nychka and Saltzman ([1998]) described the choice of knots for multivariate smoothing and kriging. In the examples we use $\mathbf{1}_d$ to denote a $d \times 1$ vector of ones.

EXAMPLE 1.   Random intercept:

$$(\mathbf{X}\boldsymbol{\beta} + \mathbf{Zu})_{ij} = \beta_0 + U_i + \beta_1 x_{ij},$$
$$1 \le j \le n_i, \ 1 \le i \le m,$$
$$\mathbf{X}_i^R = \mathbf{1}_{n_i}, \qquad \mathbf{X}^G = [x_{ij}],$$
$$\mathbf{Z}^G = \mathbf{Z}^C = \varnothing, \quad \boldsymbol{\Sigma}^R = \sigma_u^2.$$

EXAMPLE 2.   Random intercept and slope:

$$(\mathbf{X}\boldsymbol{\beta} + \mathbf{Zu})_{ij} = \beta_0 + U_i + (\beta_1 + V_i)x_{ij},$$
$$1 \le j \le n_i, \ 1 \le i \le m,$$
$$\mathbf{X}_i^R = \begin{bmatrix} 1 & x_{i1} \\ \vdots & \vdots \\ 1 & x_{in_i} \end{bmatrix}, \quad \mathbf{X}^G = \mathbf{Z}^G = \mathbf{Z}^C = \varnothing,$$
$$\boldsymbol{\Sigma}^R = \begin{bmatrix} \sigma_u^2 & \rho_{uv}\sigma_u\sigma_v \\ \rho_{uv}\sigma_u\sigma_v & \sigma_v^2 \end{bmatrix}.$$

EXAMPLE 3.   Crossed random effects model:

$$(\mathbf{X}\boldsymbol{\beta} + \mathbf{Zu})_{ii'} = \beta_0 + U_i + U_{i'}',$$
$$1 \le i \le n, \ 1 \le i' \le n',$$
$$\mathbf{X}^G = \mathbf{1}_{nn'}, \quad \mathbf{Z}^G = [\mathbf{I}_n \otimes \mathbf{1}_{n'} | \mathbf{1}_n \otimes \mathbf{I}_{n'}],$$
$$\mathbf{X}^R = \mathbf{Z}^R = \mathbf{Z}^C = \varnothing,$$
$$\mathbf{u}^G = [U_1, \ldots, U_n, U_1', \ldots, U_{n'}']^\top,$$
$$\mathrm{Cov}(\mathbf{u}^G) = \mathrm{blockdiag}(\sigma_u^2 \mathbf{I}_n, \sigma_{u'}^2 \mathbf{I}_{n'}).$$

An alternative representation of this model is

$$\mathbf{X}_i^R = \mathbf{1}_{n' \times 1}, \quad \mathbf{Z}^G = [\mathbf{1}_n \otimes \mathbf{I}_{n'}], \quad \mathbf{X}^G = \mathbf{Z}^C = \varnothing,$$
$$\boldsymbol{\Sigma}^R = \sigma_u^2, \quad \mathrm{Cov}(\mathbf{u}^G) = \sigma_{u'}^2 \mathbf{I}_{n'}.$$

This allows for implementation of hierarchical centering as described in Section [2.3].

EXAMPLE 4.   Nested random effects model:

$$(\mathbf{X}\boldsymbol{\beta} + \mathbf{Zu})_{ijk} = \beta_0 + U_i + V_{j(i)} + \beta_1 x_{ijk},$$
$$1 \le i \le m, \ 1 \le j \le n, \ 1 \le k \le p,$$
$$\mathbf{X}^G = [1 \quad x_{ijk}]_{1 \le i \le m, 1 \le j \le n, 1 \le k \le p},$$

$$\mathbf{Z}^G = [\mathbf{I}_m \otimes \mathbf{1}_{np} | \mathbf{I}_m \otimes (\mathbf{I}_n \otimes \mathbf{1}_p)],$$
$$\mathbf{X}^R = \mathbf{Z}^R = \mathbf{Z}^C = \varnothing,$$
$$\mathbf{u}^G = [U_1, \ldots, U_m, V_{1(1)}, \ldots,$$
$$V_{n(1)}, \ldots, V_{1(m)}, \ldots, V_{n(m)}]^\top,$$
$$\mathrm{Cov}(\mathbf{u}^G) = \mathrm{blockdiag}(\sigma_u^2 \mathbf{I}_m, \sigma_v^2 \mathbf{I}_{np}).$$

EXAMPLE 5.   Generalized scatterplot smoothing:

$$(\mathbf{X}\boldsymbol{\beta} + \mathbf{Zu})_i = \beta_0 + \beta_1 x_i + \sum_{k=1}^K u_k(x_i - \kappa_k)_+,$$
$$\mathbf{X}^G = [1 \quad x_i]_{1 \le i \le n},$$
$$\mathbf{Z}^G = \left[ (x_i - \kappa_k)_{+ \atop 1 \le k \le K} \right]_{1 \le i \le n},$$
$$\mathbf{X}^R = \mathbf{Z}^R = \mathbf{Z}^C = \varnothing,$$
$$\mathrm{Cov}(\mathbf{u}^G) = \sigma_u^2 \mathbf{I}_K.$$

EXAMPLE 6.   Generalized additive model:

$$(\mathbf{X}\boldsymbol{\beta} + \mathbf{Zu})_i = \beta_0 + \beta_s s_i + \sum_{k=1}^{K^s} u_k^s(s_i - \kappa_k^s)_+$$
$$+ \beta_t t_i + \sum_{k=1}^{K^t} u_k^t(t_i - \kappa_k^t)_+,$$
$$\mathbf{X}^G = [1 \quad s_i \quad t_i]_{1 \le i \le n},$$
$$\mathbf{Z}^G = \left[ (s_i - \kappa_k^s)_{+ \atop 1 \le k \le K^s} \quad (t_i - \kappa_k^t)_{+ \atop 1 \le k \le K^t} \right]_{1 \le i \le n},$$
$$\mathbf{X}^R = \mathbf{Z}^R = \mathbf{Z}^C = \varnothing,$$
$$\mathrm{Cov}(\mathbf{u}^G) = \mathrm{blockdiag}(\sigma_{us}^2 \mathbf{I}_{K^s}, \sigma_{ut}^2 \mathbf{I}_{K^t}).$$

EXAMPLE 7.   Generalized additive semiparametric mixed model:

$$(\mathbf{X}\boldsymbol{\beta} + \mathbf{Zu})_{ij} = \beta_0 + U_i + (\beta_q + V_i)q_{ij}$$
$$+ (\beta_r + W_i)r_{ij} + \beta_1 x_{ij}$$
$$+ \beta_s s_{ij} + \sum_{k=1}^{K^s} u_k^s(s_{ij} - \kappa_k^s)_+$$
$$+ \beta_t t_{ij} + \sum_{k=1}^{K^t} u_k^t(t_{ij} - \kappa_k^t)_+,$$
$$\mathbf{X}_i^R = \begin{bmatrix} 1 & q_{i1} & r_{i1} \\ \vdots & \vdots & \vdots \\ 1 & q_{in_i} & r_{in_i} \end{bmatrix},$$
$$\mathbf{X}^G = [s_{ij} \quad t_{ij} \quad x_{ij}]_{1 \le j \le n_i, 1 \le i \le m},$$

$$\mathbf{Z}^G = \begin{bmatrix} (s_{ij} - \kappa_k^s)_+ & (t_{ij} - \kappa_k^t)_+ \\ {}_{1 \le k \le K^s} & {}_{1 \le k \le K^t} \end{bmatrix},$$

$$\mathbf{Z}^C = \varnothing,$$

$$\mathbf{\Sigma}^R = \text{unstructured } 3 \times 3 \text{ covariance matrix},$$

$$\text{Cov}(\mathbf{u}^G) = \text{blockdiag}(\sigma_{us}^2 \mathbf{I}_{K^s}, \sigma_{ut}^2 \mathbf{I}_{K^t}).$$

EXAMPLE 8. Generalized bivariate smoothing/low-rank kriging:

$$(\mathbf{X}\boldsymbol{\beta} + \mathbf{Z}\mathbf{u})_i = \beta_0 + \boldsymbol{\beta}_1^\top \mathbf{x}_i + \sum_{k=1}^K u_k C(\|\mathbf{x}_i - \boldsymbol{\kappa}_k\|),$$

$$\mathbf{X}^G = [\, 1 \quad \mathbf{x}_i^\top \,]_{1 \le i \le n},$$

$$\mathbf{Z}^G = \left[ C\left( \|\mathbf{x}_i - \boldsymbol{\kappa}_k\| \right) \right]_{\substack{1 \le k \le K \\ 1 \le i \le n}},$$

$$\mathbf{X}^R = \mathbf{Z}^R = \mathbf{Z}^C = \varnothing,$$

$$\text{Cov}(\mathbf{u}^G) = \sigma_u^2 \mathbf{I}.$$

Here $\|\mathbf{v}\| \equiv \sqrt{\mathbf{v}^\top \mathbf{v}}$, $C(r) = r^2 \log |r|$ corresponds to low-rank thin plate splines with smoothness parameter set to 2 (as defined in Wahba, 1990) and $C(r) = \exp(-|r/\rho|)(1 + |r/\rho|)$ corresponds to Matérn low-rank kriging with range $\rho > 0$ and smoothness parameter set to $3/2$ (as defined in Stein, 1999; Kammann and Wand, 2003). Several more examples could be added, including some where $\mathbf{Z}^C \ne \varnothing$.

## 2.1 Fixed Effects Priors

Throughout we take the prior distribution of the fixed effects vector $\boldsymbol{\beta}$ to be of the form

$$[\boldsymbol{\beta}] \sim N(\mathbf{0}, \mathbf{F})$$

for some covariance matrix $\mathbf{F}$. In practice it is common to take $\mathbf{F}$ to be diagonal with very large entries, corresponding to noninformative priors on the entries of $\boldsymbol{\beta}$. Such a strategy ensures that, with appropriate choice of prior for the variance components, the resulting joint posterior distribution of the parameters will be proper. Even so, the resulting posterior distributions approximate those based on uniform priors for $\boldsymbol{\beta}$. For normal models, Gelman (2005) noted that because we typically have enough data to estimate these coefficients from the data, any noninformative prior is adequate. For binary response models, Natarajan and Kass (2000) showed that, under mild regularity conditions that usually amount to soft requirements on the number of successes and failures in the data set, use of a uniform distribution for $\boldsymbol{\beta}$ in conjunction with an appropriate prior for the variance components results in a proper posterior. For logistic regression, Bedrick, Christensen and Johnson (1997) noted that the normal prior for $\boldsymbol{\beta}$ is convenient in large sample situations in which the posterior is approximately normal. In other situations, one should be cautious about using a normal prior with large covariances, because the induced prior distributions for each $P(y = 1)$ can have point masses at zero and one. In such cases, it may be preferable to use the conditional means priors proposed by Bedrick, Christensen and Johnson (1996), which specify prior distributions on the success probabilities directly.

## 2.2 Covariance Matrix Priors

Over the last decade-and-a-half, prior elicitation for the variance components in Bayesian GLMMs has been an active area of statistical research. Several authors have demonstrated that the use of improper priors for these parameters can lead to improper posteriors, with Gibbs samplers unable to detect such ill-conditioning (Hobert and Casella, 1996). As a result, a popular choice is the use of proper but "diffuse" conditionally conjugate priors. In the GLMM setting with normal random effects, this corresponds to an inverse gamma (IG) distribution for a single variance component and an inverse Wishart distribution for a variance–covariance matrix. For hierarchical versions of GLMMs, however, recent research has shown that these priors can actually be quite informative, leading to inferences that are sensitive to choice of the hyperparameters for these distributions (Natarajan and McCulloch, 1998; Natarajan and Kass, 2000; Gelman, 2005). Natarajan and Kass (2000) and Gelman (2005) have proposed alternative prior elicitation strategies that improve upon the conditionally conjugate priors. In Section 4 we outline a sensitivity analysis approach that takes these latest proposals into account.

## 2.3 Hierarchical Centering

Hierarchical centering of parameters can improve convergence of Markov chain Monte Carlo schemes (Section 3) for fitting Bayesian mixed models (e.g., Gelfand, Sahu and Carlin, 1995). In the context of this section, hierarchical centering involves reparametrization of $(\boldsymbol{\beta}^R, \mathbf{u}^R)$ to $(\boldsymbol{\beta}^R, \boldsymbol{\gamma})$, where

$$\boldsymbol{\gamma} \equiv \{(\mathbf{Z}^R)^\top \mathbf{Z}^R\}^{-1} (\mathbf{Z}^R)^\top \mathbf{X}^R \boldsymbol{\beta}^R + \mathbf{u}^R.$$

The new vector of parameters $\boldsymbol{\gamma}$ can be further divided into $m$ subvectors $\boldsymbol{\gamma}_i$ with $\boldsymbol{\gamma}_i = \boldsymbol{\beta}^R + \mathbf{u}_i^R$, so

that

$$\boldsymbol{\gamma} = \begin{bmatrix} \gamma_1 \\ \vdots \\ \gamma_m \end{bmatrix}.$$

Then the general design generalized linear mixed model becomes

$$\mathbf{X}\boldsymbol{\beta} + \mathbf{Z}\mathbf{u} = \mathbf{Z}^R\boldsymbol{\gamma} + \mathbf{X}^G\boldsymbol{\beta}^G + \sum_{\ell=1}^{L} \mathbf{Z}_\ell^G \mathbf{u}_\ell^G + \mathbf{Z}^C \mathbf{u}^C.$$

Note that hierarchical centering is not a well-defined concept for general design or spatial structures, because $\mathbf{u}^G$ and $\mathbf{u}^C$ cannot be centered in a hierarchical way similarly to that for $\mathbf{u}^R$. As a result, the general design and spatial structures do not contribute to the model for the mean in a conditionally hierarchical manner.

## 2.4 Applications

This section describes three public health applications that benefit from general design Bayesian GLMM analysis. The analyses are postponed to Section 4.

2.4.1 *Respiratory infection in Indonesian children.* Our first example involves longitudinal measurements on 275 Indonesian children. Analyses of these data have appeared previously in the literature (e.g., Diggle, Liang and Zeger, 1994; Lin and Carroll, 2001), so our description of them will be brief. The response variable is binary: the indicator of respiratory infection. The covariate of most interest is the indicator of vitamin A deficiency. However, the age of the child has been seen in previous analyses to have a nonlinear effect.

A plausible model for these data is the Bayesian logistic additive mixed model

$$(5) \quad \begin{aligned} \text{logit}\{P(\texttt{respiratory infection}_{ij} = 1)\} \\ = \beta_0 + U_i + \boldsymbol{\beta}^\top \mathbf{x}_{ij} + f(\texttt{age}_{ij}), \end{aligned}$$

where $1 \le i \le 275$ indexes child and $1 \le j \le n_i$ indexes the repeated measures within child. Here $U_i \overset{\text{ind.}}{\sim} N(0, \sigma_U^2)$ is a random child effect, $\mathbf{x}_{ij}$ denote measurements on a vector of nine covariates—height and indicators for vitamin A deficiency, sex, stunting and visit number—and $f$ is modeled using penalized splines with spline basis coefficients $u_k \overset{\text{ind.}}{\sim} N(0, \sigma_u^2)$.

2.4.2 *Caregiver stress and respiratory health.* The Home Allergen and Asthma study is an ongoing longitudinal study that is investigating risk factors for incidence of childhood respiratory problems including asthma, allergy and wheeze (Gold et al., 1999). The portion of the study data that we will consider consists of 483 families who were followed for two and-a-half years after the birth of a child. At the start of the study, a number of demographic variables were measured on each family, including race, categorized household income, categorized caregiver educational level and child's gender. Additionally, one of the hypothesized risk factors for childhood respiratory problems is exposure to a stressful environment (Wright et al., 2004). Each child's environmental stress level was measured approximately bimonthly by a telephone interview and assessed on a discrete ordinal scale from 0 (no stress) to 16 (very high stress). This assessment was based on the four-item Perceived Stress Scale (PSS-4) (Cohen, 1988).

Let $1 \le i \le 483$ index family and let $1 \le j \le n_i$ index the repeated measurements within each family. We arrived at the following Bayesian Poisson additive mixed model for stress experience by caregiver $i$ when the child was $\texttt{age}_{ij}$:

$$(6) \quad \begin{aligned} \texttt{stress}_{ij} \\ \sim \text{Poisson}[\exp\{\beta_0 + U_i + \boldsymbol{\beta}^\top \mathbf{x}_{ij} + f(\texttt{age}_{ij})\}]. \end{aligned}$$

The random intercept, $U_i \overset{\text{ind.}}{\sim} N(0, \sigma_U^2)$, is a random family effect, and $\mathbf{x}_i$ includes indicators of annual family income and race (see Figure 4 for details). The term $f(\texttt{age}_{ij})$ is a nonparametric term that we model using penalized splines with spline coefficients $u_k \overset{\text{ind.}}{\sim} N(0, \sigma_u^2)$. We include the nonparametric term in the model for the effect of stress as a function of child's age because, outside of anecdotal evidence, we do not know of a biologically motivated parametric model for stress as a function of child's age. We arrived at the other terms in the model (and removed other demographic terms and interactions from the model) based on discussions with the investigators in the study and exploratory data analyses that we fitted via maximum PQL.

2.4.3 *Standardized cancer incidence and proximity to a pollution source.* Elevated cancer rates were observed in a region of Massachusetts, USA, known as Upper Cape Cod, during the mid-1980s, and one risk factor of interest is a fuel dump at the Massachusetts Military Reservation (MMR) (Kammann and Wand, 2003; French and Wand, 2004). For nearly 20 years the Massachusetts Department of Public Health (MDPH) has maintained a cancer registry data base which records incident cases for 22 types of cancers, including lung, breast and prostate cancers. In this example we focus on female lung cancer between 1986 and 1994.

We use a semiparametric Poisson spatial model to investigate the relationship between census tract level female lung cancer standardized incidence rates (SIRs) and distance to the MMR. Let $i = 1, \ldots, 45$ represent the census tracts in the study, and let observed$_i$ and expected$_i$ be the observed and expected number of incident cases of female lung cancer in tract $i$ (i.e., numerator and denominator of the SIR), respectively. After fitting a number of models that included terms for additional demographic factors and water source, we arrived at the semiparametric Poisson spatial model

$$
\begin{aligned}
\text{observed}_i \\
(7) \quad &\sim \text{Poisson}[\text{expected}_i \exp\{\beta_0 + U_i^C + \beta_1 x_i \\
&\qquad + f(\text{dist}_i)\}],
\end{aligned}
$$

where $x_i$ is the percentage of women in tract $i$ who were over 15 and employed outside the home in 1989, and dist$_i$ is the distance from the centroid of census tract $i$ to the centroid of the MMR. Here, $\mathbf{u}^C = (U_1^C, \ldots, U_{45}^C)^\top$ is a vector of spatially correlated random effects with intrinsic Gaussian autoregression distribution parametrized by variance component $\sigma_c^2$, as defined in (4). To complete the specification of the spatial correlation model, we choose a cutoff distance value $d$ and treat two census tracts as neighbors if the distance between their centroids is less than or equal to $d$. We choose $d = 7.5$ km, which corresponds to the cutoff such that every census tract has at least one neighbor. We model the nonparametric term $f(\text{dist}_i)$ using penalized splines with coefficients $u_k \stackrel{\text{ind.}}{\sim} N(0, \sigma_u^2)$.

## 3. FITTING VIA MARKOV CHAIN MONTE CARLO

In the general design GLMM (1) and (2), the posterior distribution of

$$
\boldsymbol{\nu}^\top \equiv [\boldsymbol{\beta}^\top \quad \mathbf{u}^\top]
$$

is

$$
[\boldsymbol{\nu}|\mathbf{y}] = \left( \int \exp\{\mathbf{y}^\top \mathbf{C}\boldsymbol{\nu} - \mathbf{1}^\top b(\mathbf{C}\boldsymbol{\nu})\right.
$$

$$
(8) \qquad \left. - \tfrac{1}{2}(\log|\mathbf{G}| + \boldsymbol{\nu}^\top \mathbf{V}^{-1}\boldsymbol{\nu})\} [\mathbf{G}] \, d\mathbf{G}\right)
$$

$$
\cdot \left( \int\int \exp\{\mathbf{y}^\top \mathbf{C}\boldsymbol{\nu} - \mathbf{1}^\top b(\mathbf{C}\boldsymbol{\nu})\right.
$$

$$
- \tfrac{1}{2}(\log|\mathbf{G}| + \boldsymbol{\nu}^\top \mathbf{V}^{-1}\boldsymbol{\nu})\}
$$

$$
\left. \cdot [\mathbf{G}] \, d\mathbf{G} \, d\boldsymbol{\nu}\right)^{-1},
$$

where $\mathbf{C} \equiv [\mathbf{X} \quad \mathbf{Z}]$, $\mathbf{V} \equiv \text{blockdiag}(\mathbf{F}, \mathbf{G})$ and $[\mathbf{G}]$ is the prior on the variance components in $\mathbf{G}$. These integrals are analytically intractable for most problems. Furthermore, in the applications we consider, the dimensionality precludes the use of numerical integration. A standard remedy is to apply a Markov chain Monte Carlo algorithm to draw samples from (8). An overview of MCMC is provided by Gilks, Richardson and Spiegelhalter (1996).

The MCMC methods break up the model parameters into subsets and then sample from the conditional distributions given the remaining parameters and data, which are often called full conditionals. In the general design GLMM, the natural breakdown of the parameters is into $\boldsymbol{\nu}$ and $\mathbf{G}$, leading to the full conditionals

$$
[\boldsymbol{\nu}|\mathbf{G}, \mathbf{y}] \quad \text{and} \quad [\mathbf{G}|\boldsymbol{\nu}, \mathbf{y}].
$$

The latter full conditional has a standard form when the prior on the variance components is inverse gamma or Wishart, which are "conditionally conjugate" priors for this model, but not when, say, a folded-Cauchy prior is used (e.g., Gelman, 2005). The first full conditional has the general form

$$
[\boldsymbol{\nu}|\mathbf{G}, \mathbf{y}] \propto \exp\{\mathbf{y}^\top \mathbf{C}\boldsymbol{\nu} - \mathbf{1}^\top b(\mathbf{C}\boldsymbol{\nu}) - \tfrac{1}{2}\boldsymbol{\nu}^\top \mathbf{V}^{-1}\boldsymbol{\nu}\},
$$

which is a nonstandard distribution unless $\mathbf{y}$ is conditionally Gaussian. Clever strategies such as adaptive rejection sampling (Gilks and Wild, 1992) and slice sampling (e.g., Besag and Green, 1993; Neal, 2003) are required to draw samples. The most common versions of these algorithms work with the full conditionals of the components $\boldsymbol{\nu}$. When $\mathbf{V}$ is diagonal, these full conditionals are of the form

$$
\begin{aligned}
[\nu_k|\boldsymbol{\nu}_{-k}, \mathbf{G}, \mathbf{y}] \\
(9) \quad &\propto \exp\{(\mathbf{C}^\top \mathbf{y})_k \nu_k - \mathbf{1}^\top b(\mathbf{c}_k \nu_k + \mathbf{C}_{-k}\boldsymbol{\nu}_{-k}) \\
&\qquad - \tfrac{1}{2}\nu_k^2/(\mathbf{V})_{kk}\}.
\end{aligned}
$$

Here $\mathbf{c}_k$ is the $k$th column of $\mathbf{C}$, $\mathbf{C}_{-k}$ is $\mathbf{C}$ with the $k$th column omitted, $\nu_k$ is the $k$th entry of $\boldsymbol{\nu}$ and $\boldsymbol{\nu}_{-k}$ is $\boldsymbol{\nu}$ with the $k$th entry omitted. It is easily shown that (9) is log-concave, which permits use of adaptive rejection sampling and simplifies slice sampling. These algorithms can also be used to sample from the full conditionals for the variance components when necessary.

Zhao (2003) provides a detailed account of MCMC for general design GLMM and compares several strategies via simulation. One of the conclusions drawn from the simulations is that the `WinBUGS` package (Spiegelhalter, Thomas and Best, 2000) performs excellently among various "off-the-shelf" competitors. This is very good news because it saves the user having to write his or her own MCMC code. However, it should be noted that for large models `WinBUGS` can take quite some time to obtain a fit. Also, the analysis must be performed on a particular platform (`Windows`). Assuming that computation time is not an issue and that `Windows` is available, we can report that fitting of general design GLMMs via `WinBUGS` has a large chance of success. For our analyses we had access to several personal computers and ran multiple chains in parallel to assess convergence and prior sensitivity. This reduced the elapsed time it took to compute each one of the analyses by an order of magnitude.

## 4. DATA ANALYSES

### 4.1 Input Values and Prior Distributions

We used `WinBUGS` to fit the models described in Section 2.4. However, several input values and prior distributions needed to be specified, so we preface the analyses with the particular choices that were made. A more detailed study on the use of `WinBUGS` for fitting models of this type is provided by Crainiceanu, Ruppert and Wand (2005).

Based on the recommendations of Gelfand, Sahu and Carlin (1995), we used hierarchical centering of random effects. All continuous covariates were standardized to have zero mean and unit standard deviation. A strategy such as this is necessary for the method to be scale invariant given fixed choices for the hyperparameters. We used radial cubic basis functions for smooth function components. Apart from making the fitted functions smooth and requiring a relatively small number of knots, Crainiceanu, Ruppert and Wand (2005) reported that they had good mixing properties in MCMC analysis. Radial cubic basis function modeling of a function $f$ entails putting $f(x) = \beta_0 + \beta_1 x + \mathbf{Z}_x \mathbf{u}$, where

$$\mathbf{Z}_x = \left[|x - \kappa_k|^3\right]_{1 \le k \le K} \left[|\kappa_k - \kappa_{k'}|^3\right]_{1 \le k,k' \le K}^{-1/2} \quad \text{and}$$

$$(10) \quad \mathbf{u} \sim N(0, \sigma_u^2 \mathbf{I})$$

(French, Kammann and Wand, 2001), with $\kappa_k = (\frac{k+1}{K+2})$th quantile of the unique predictor values. In general, $K$ can be chosen using rules such as

$$K = \min(\tfrac{1}{4}(\text{number of unique predictor values}), 35)$$

or those given in Ruppert (2002). However, often considerably smaller $K$ can be used through the benefit of faster MCMC fitting. This approach was taken in our analyses.

We considered several common variance component priors. These were inverse gamma with equal scale and shape,

$$[\sigma^2] \propto (\sigma^2)^{-(a+1)} e^{-a/\sigma^2},$$

denoted by $\mathrm{IG}(a, a)$ for $a = 0.001, 0.01, 0.1$ (Spiegelhalter et al., 2003), and the folded-$t$ class of priors for $\sigma$ (Gelman, 2005)

$$[\sigma] \propto \left(1 + \frac{1}{\nu}\left(\frac{\sigma}{s}\right)^2\right)^{-(\nu+1)/2},$$

where $s$ and $\nu$ are fixed scale and degrees of freedom hyperparameters, respectively. We investigated the sensitivity of the model fit to the choice of the hyperparameters. Results showed that fits based on the IG priors were stable for $a \ge 0.01$, but those obtained assuming $a = 0.001$ behaved erratically. Out of the remaining choices, the folded-Cauchy prior, a member of the folded-$t$ class of priors, performed well.

As a result of these empirical comparisons, in our examples we take the approach of fitting models under multiple prior distributions for the variance components and assessing the sensitivity of the results to these assumptions. Due to its popularity, we fit general design GLMMs using independent IG(0.01, 0.01) priors for each variance component. Results suggest this prior performs well for the examples considered in this paper. We also refitted the models using independent folded-Cauchy prior distributions for each variance component. For a variance component square root $\sigma$ and fixed scale parameter $s$, the folded-Cauchy distribution has probability distribution

$$[\sigma] \propto (\sigma^2 + s^2)^{-1}.$$

Following Gelman (2005), we take $s = 25$ in our examples and check the sensitivity of results to this choice by also fitting the models for $s = 12$. This prior can be implemented in `WinBUGS` using the flexible feature that allows a user to code an arbitrary prior distribution for the model parameters (see the Appendix). We also ran the models using a Uniform(0, 100) prior on $\sigma$. A theoretical comparison of such priors in the general design GLMM setting is a topic worthy of future research.

Table 1 summarizes the input values and prior distributions that were used.

## 4.2 Respiratory Infection in Indonesian Children

Using the prior distributions and input values given in Table 1, `WinBUGS` produced the output for the $\boldsymbol{\beta}$ coefficients summarized in Figure 1. It is seen that, for this model, the chains mix quite well with little significant autocorrelation and Gelman–Rubin $\sqrt{R}$ values (Gelman and Rubin, 1992) all less than 1.01. Vitamin A deficiency is seen to have a borderline positive effect on respiratory infection, which is in keeping with previous analyses. Similar comments apply to sex and some of the visit numbers.

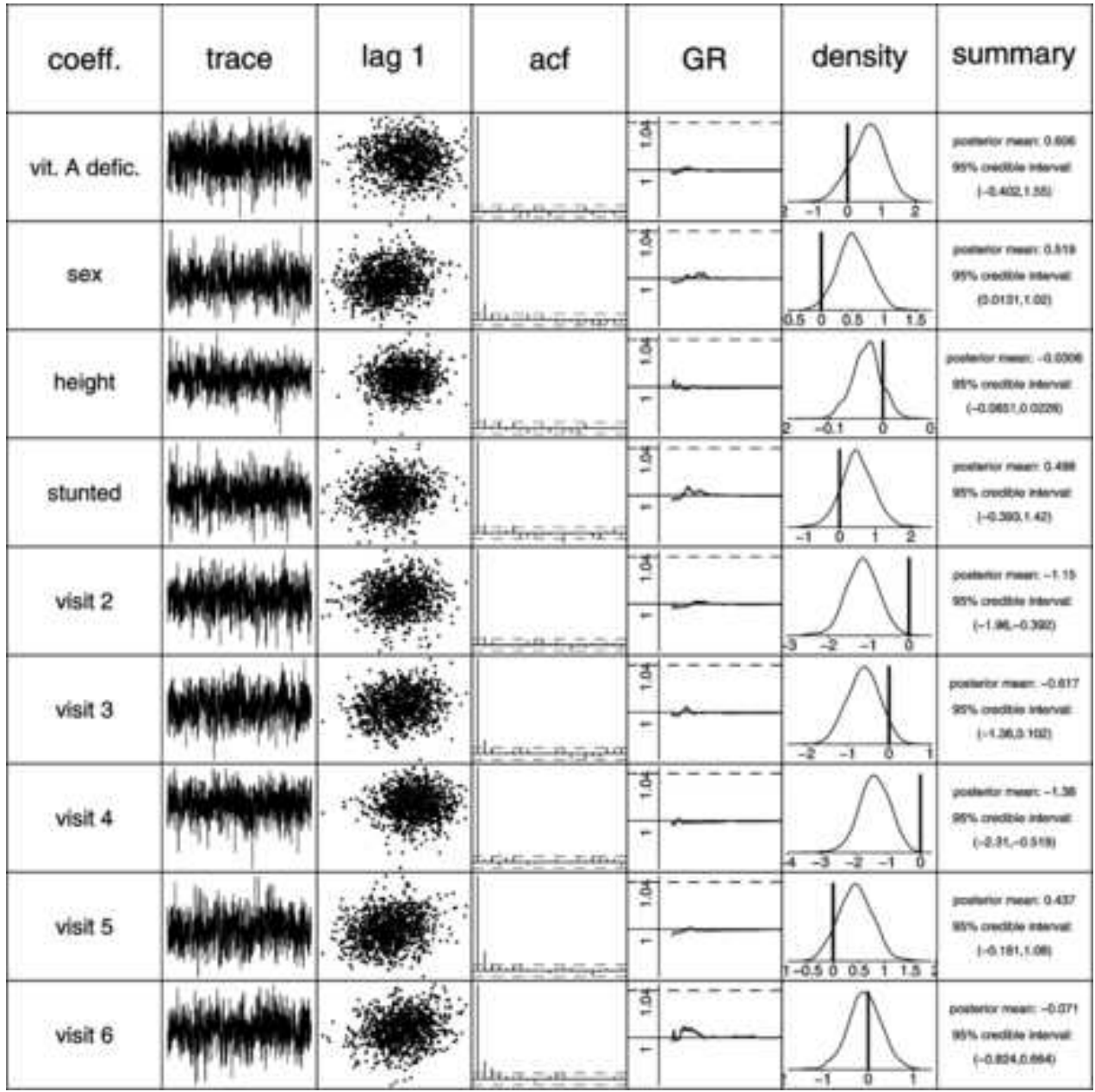

FIG. 1. *Summary of `WinBUGS` output for parametric components of* (5). *The full titles for columns are name of variable, trace plot of sample of corresponding coefficient, plot of sample against 1-lagged sample, sample autocorrelation function, Gelman–Rubin $\sqrt{\widehat{R}}$ diagnostic, kernel estimate of posterior density and basic numerical summaries.*

TABLE 1
*Input values and prior distributions used in `WinBUGS` for analyses*

| | |
|---|---|
| Hierarchical centering used for random intercepts. | |
| Continuous covariates standardized to have | |
| zero mean and unit standard deviation. | |
| Radial cubic basis functions for smooth functions. | |
|    length of burn-in | 5000 |
|    length of "kept" chain | 5000 |
|    thinning factor | 5 |
|    prior for fixed effects | $N(0, 10^8)$ |
| prior for variance components $\begin{cases} \text{IG}(0.01, 0.01) \\ \text{folded-Cauchy with } s = 12,25 \\ \text{Uniform}(0, 100) \end{cases}$ | |

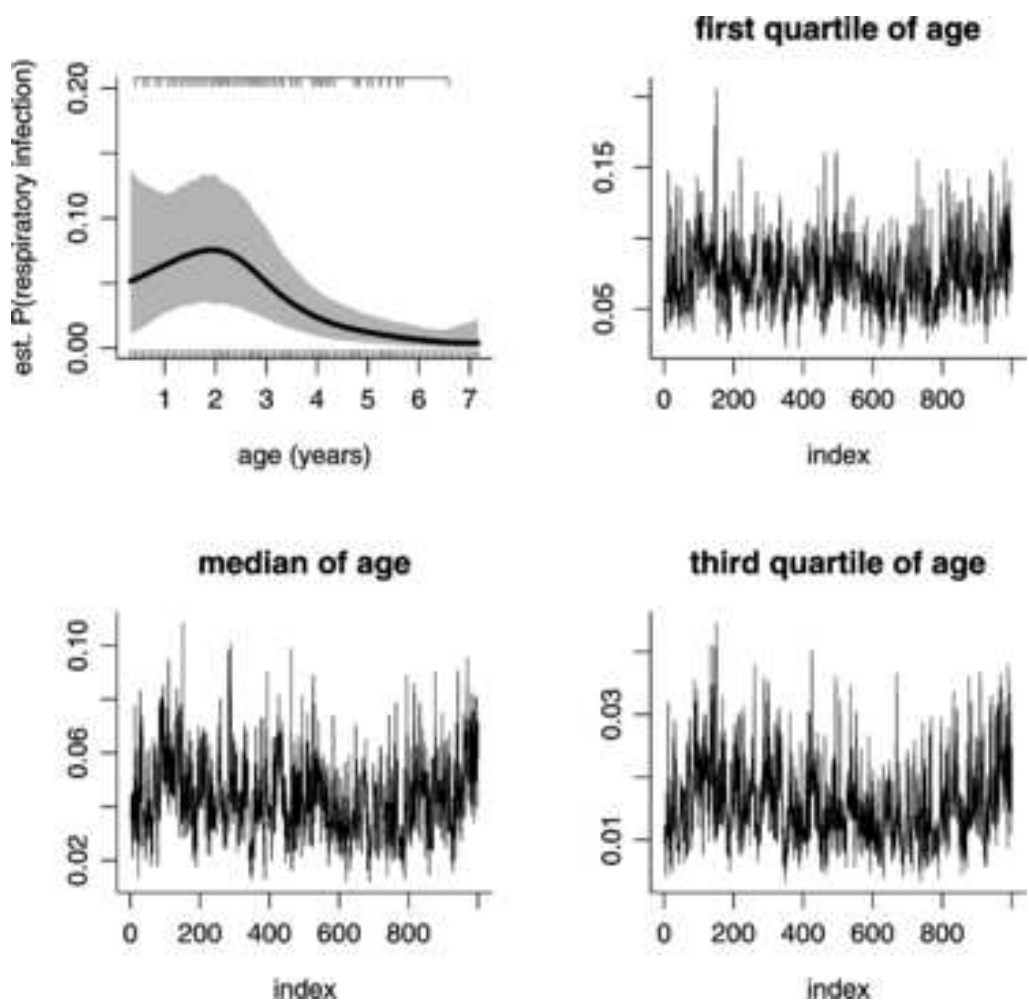

FIG. 2.  *Summary of `WinBUGS` output for estimate of `f(age)`. The top left panel is the posterior mean of the estimated probability of respiratory infection with all other covariates set to their average values. The shaded region is a corresponding pointwise 95% credible set. The remaining panels are trace plots of samples used to produce the top left plot at quartiles of the age data.*

The estimated effect of age is summarized in the top left panel of Figure 2 and is seen to be significant and nonlinear. The remaining panels show good mixing of the chains corresponding to the estimated age effect at quartiles of the age data. Gelman–Rubin $\sqrt{R}$ plots (not shown here) support convergence of these chains.

To assess the sensitivity of our conclusions to the choice of variance component priors, we also ran the Gibbs samplers assuming the folded-Cauchy and Uniform priors for the random effects standard deviations (see Section 2.2). Figure 3 shows the pos-terior estimates and 95% credible intervals for the regression coefficients of interest using the default independent IG priors, independent folded-Cauchy priors with $s = 25$, independent folded-Cauchy priors with $s = 12$ and independent $U(0, 100)$ priors. This figure shows that results are not sensitive to this choice, with the changes in the posterior means never more than 2% of that obtained from the IG specification and the credible intervals never more than 6.5% wider than their IG counterparts.

### 4.3 Caregiver Stress and Respiratory Health

For this example, we also used the priors and input values given in Table 1 and we provide the `WinBUGS` code in the Appendix. For the spline we used 12 knots that were spaced evenly on the percentiles of age. We found that the fit did not change noticeably if we used more knots and we chose a small number of knots for computational efficiency.

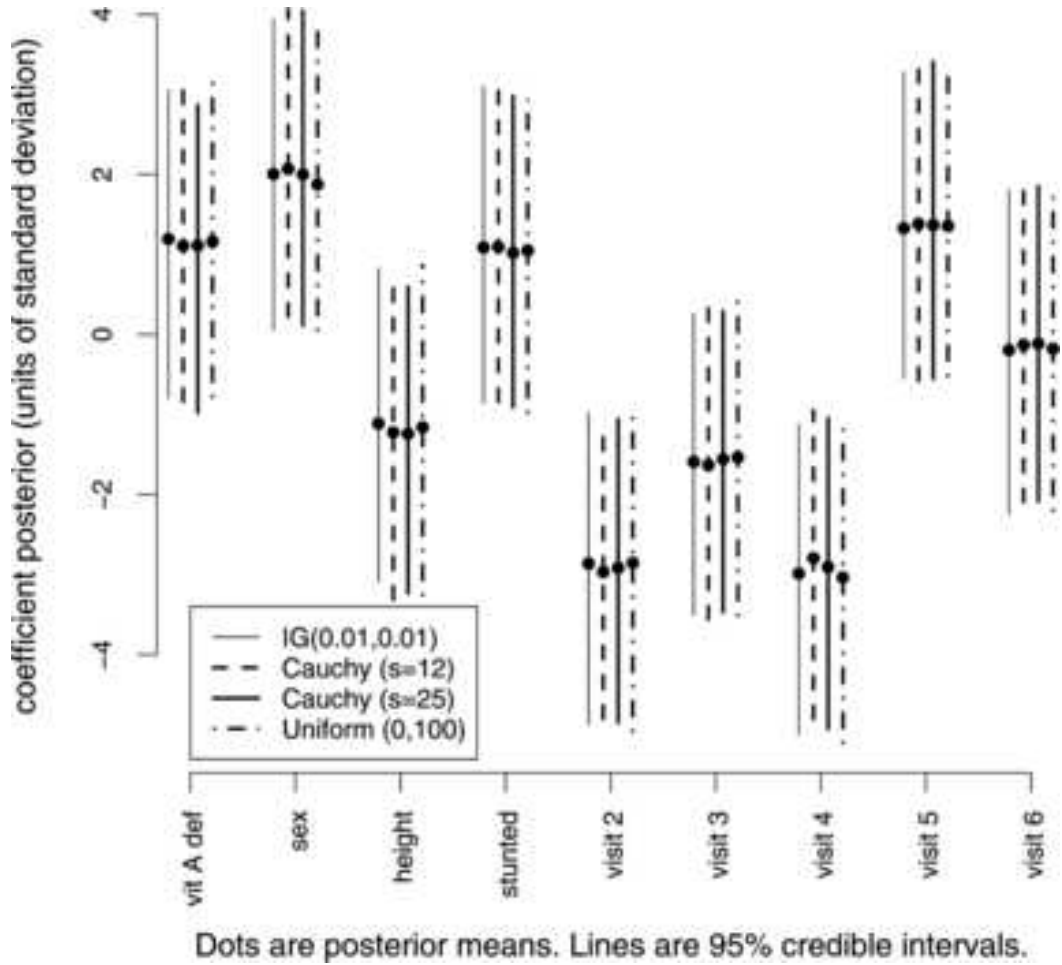

Dots are posterior means. Lines are 95% credible intervals.

FIG. 3. *Results of sensitivity analysis for variance component priors for model* (5).

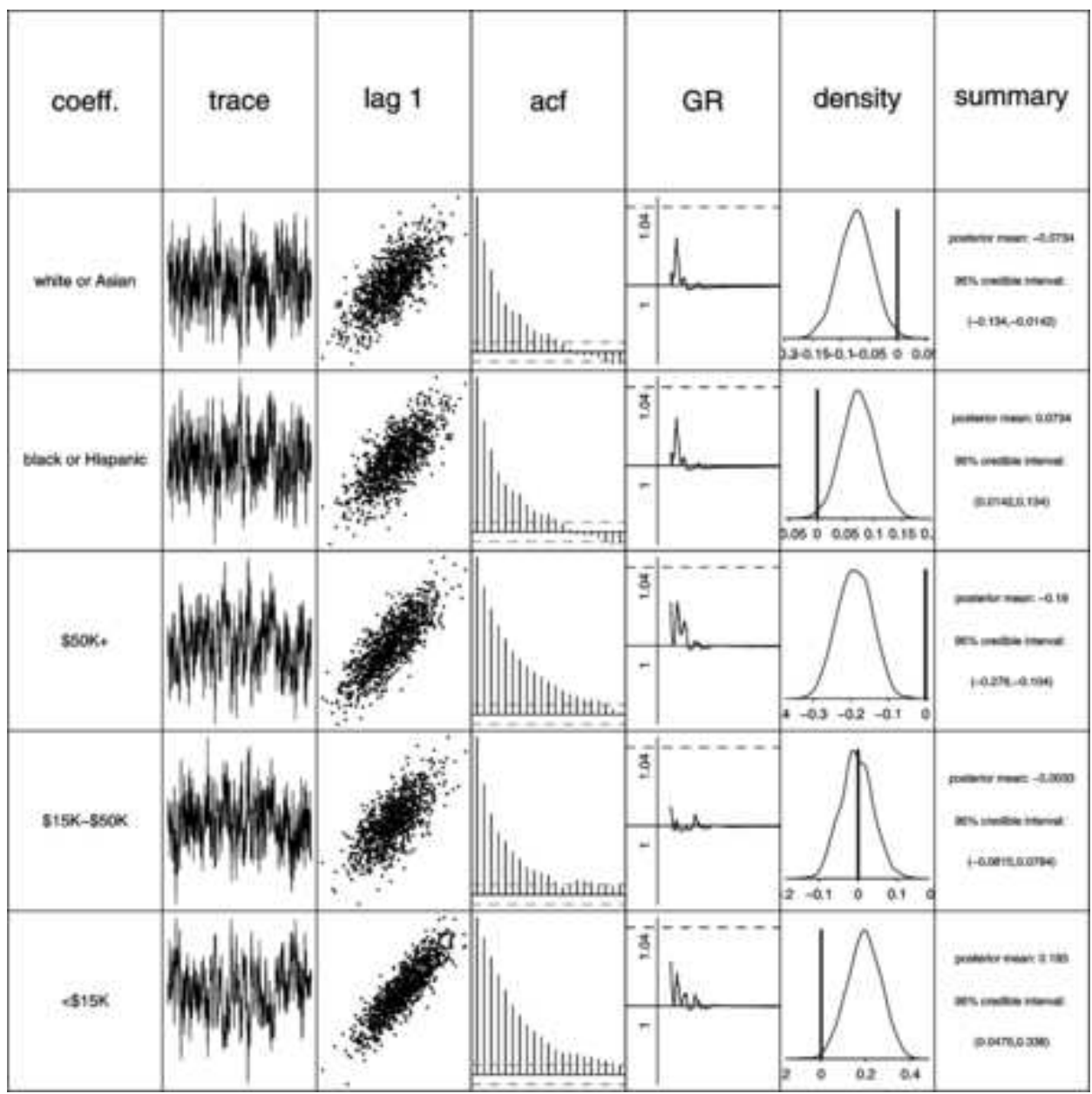

FIG. 4. *Summary of* `WinBUGS` *output for parametric components of* (6). *The full titles of columns are name of variable, trace plot of sample of corresponding coefficient, plot of sample against 1-lagged sample, sample autocorrelation function, Gelman–Rubin* $\sqrt{\hat{R}}$ *diagnostic, kernel estimate of posterior density and basic numerical summaries. The coefficients can be interpreted as time invariant offsets to the time varying population mean.*

Figure 4 shows the Bayes estimates and credible intervals for the $\boldsymbol{\beta}$ coefficients as well as an assessment of the convergence of the chains. The coefficients can be interpreted as category-specific offsets from the population mean. The chains had a moderate autocorrelation and the Gelman–Rubin $\sqrt{\hat{R}}$ values were all less than 1.04. Figure 5 contains the estimated age effect and trace plots for the effect of age at the quartiles of the data. Again, the Gelman–Rubin $\sqrt{\hat{R}}$ values were less than 1.04 and support convergence. The figures are based on the chain that used the independent inverse gamma priors for the variance components. Fits that used independent Cauchy ($s = 25$) priors for the square root of the variance components changed neither the posterior means nor the widths of the credible intervals for

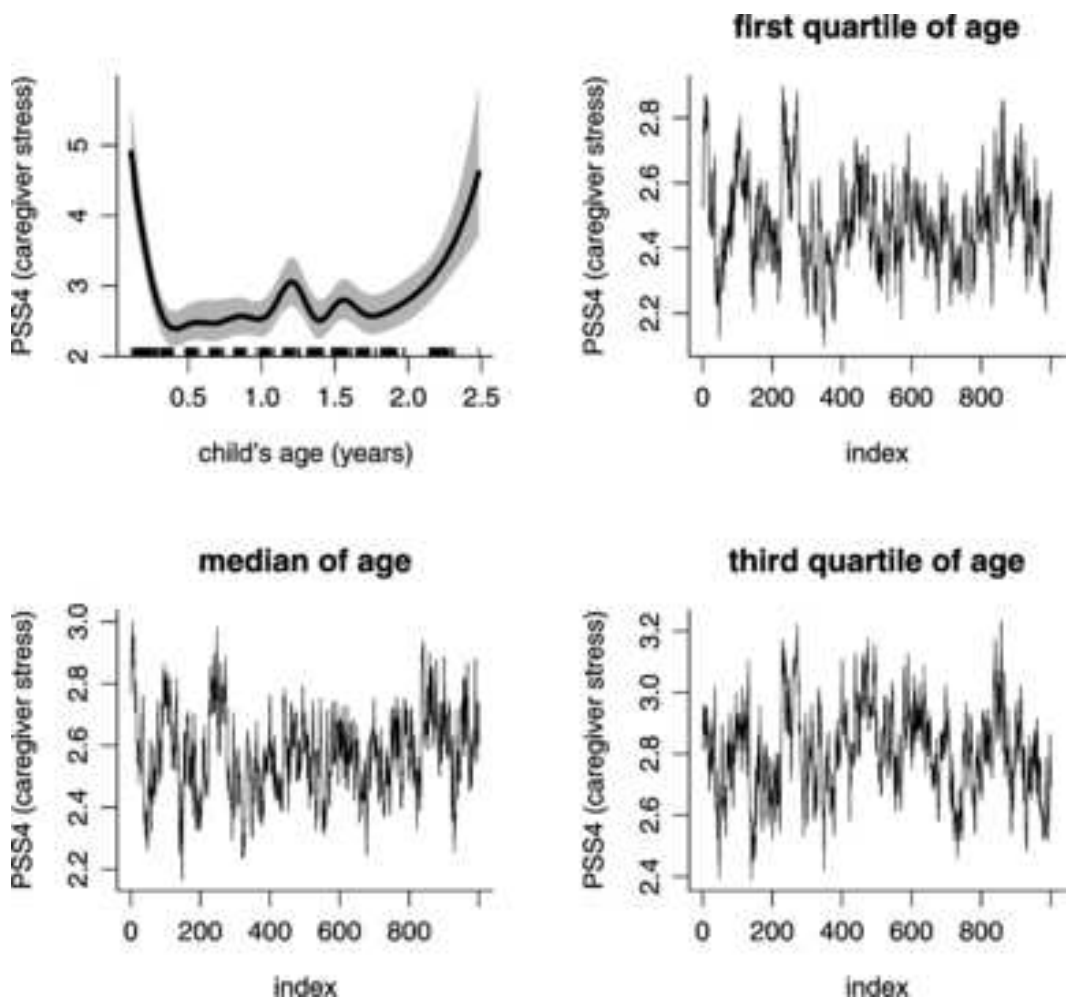

Fɪɢ. 5. *Summary of* `WinBUGS` *output for estimate of* $f(\texttt{age})$*. The top left panel is the posterior mean of the mean stress (PSS-4) as a function of age with all other covariates set to their average values. The shaded region is a corresponding pointwise 95% credible set. The remaining panels are trace plots of samples used to produce the top left plot at quartiles of the age data.*

the parameters of interest by more than 4.7%. The posterior mean and confidence set for $f(\texttt{age})$ was also relatively insensitive to the prior on the variance components in this example.

Two aspects of the fit that were of interest to the investigators in the study included the inverse dose response relationship between income and stress, and that race was significantly related to environmental stress even after accounting for the effect of income. The nonparametric estimate of stress as a function of the child's age was also interesting and suggests that relatively stressful times include the first few months, when the child is approximately a year old, and beyond age two.

### 4.4 Standardized Cancer Incidence and Proximity to a Pollution Source

As in the previous examples, we started with the prior distributions and inputs in Table 1. In this case though, the chain required a longer burn in. We found that a burn in of length 15,000 was sufficient to produce acceptable convergence. Figure 6 (bottom panel) contains the resulting convergence diagnostics and inferences for the parameters in the model. The middle panel of Figure 6 contains an estimate of the contribution of distance to the MMR to the standardized incidence and trace plots of the function estimate at the quartiles of distance. The Gelman–Rubin $\sqrt{R}$ values for the estimates at these quartiles were less than 1.04 and support convergence. Finally, the top panel of Figure 6 maps the estimated SIRs based on the model fit, demonstrating the smoothing achieved by the spatial model. The figures are based on fits that used independent IG(0.01, 0.01) priors for the variance components. Fits that used independent Cauchy ($s = 25$) priors for the square root of the variance components decreased the length of the credible interval for the effect of percent working by 6.7% and lowered the posterior mean by 3.1%. The posterior mean and confidence set for $f(\texttt{dist}_i)$ also changed very little.

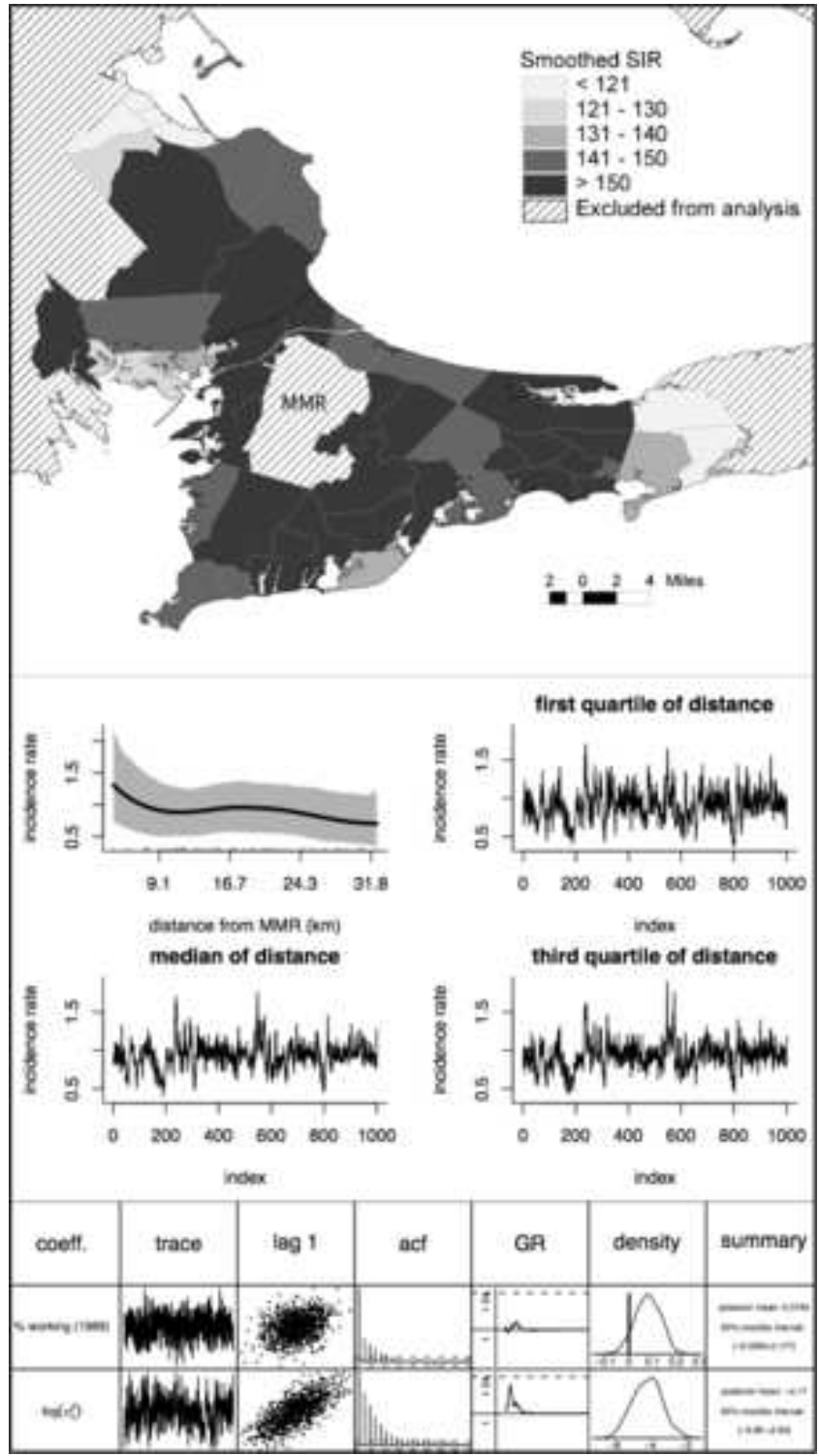

Fig. 6.  *Summary of `WinBUGS` output for the fit of (7). The top panel contains a spatial plot of the smoothed SIRs ( posterior means of the $U_i^c$'s). The middle panel shows the estimated $f(\text{dist})$ and a corresponding pointwise 95% credible set along with trace plots of the samples of the function at the quartiles of distance. The bottom panel displays summaries of other parameters of interest and convergence diagnostics. Additionally, the Gelman–Rubin $\sqrt{\hat{R}}$ diagnostics were less than 1.03 for all the $U_i^c$'s and for the $f(\text{dist})$ at the quartiles of distance. The MMR is the area in the center of the map that is excluded from the analysis.*

The fitted model suggests a nominally positive relationship between the percentage of women who were working outside the home in 1989 and standardized lung cancer incidence rates at the census tract level. Further, the estimated curve $f(\texttt{dist}_i)$ suggests an increased standardized incidence rate for census tracts that are closer than about 10 km to the MMR after controlling for other factors, and the map suggests that areas immediately east of the MMR exhibit the highest SIRs. None of the estimated effects of the model covariates is strongly significant. Regardless of statistical significance, however, we emphasize that this type of "cancer cluster" study should be viewed as exploratory since the study design is ecological (e.g., Kelsey, Whittemore, Evans and Thompson, 1996, Chapter 10). Additionally, reanalyses of similar studies have demonstrated that unmeasured confounders could radically change the conclusions in these types of analyses (e.g., Aherns et al., 2001).

## 5. DISCUSSION

As illustrated by the analyses in the previous section, general design Bayesian GLMMs are a very useful structure. In this article we have demonstrated that `WinBUGS` provides good off-the-shelf MCMC fitting of these models. Some of the reviewers have pointed out the possibility of designing MCMC algorithms that take advantage of the special structure of Bayesian GLMMs that is summarized in Section 2. We have done some exploration in this direction (Zhao, 2003), but would welcome such research from MCMC specialists. In the meantime, use of `WinBUGS` is our recommended fitting method.

## APPENDIX: WINBUGS CODE

In this Appendix we list the `WinBUGS` code used for the data analyses of Section 4. Note that the spline basis functions and hyperparameters are inputs.

The following code was used to fit (5) to the data on respiratory infection of Indonesian children. Here inverse gamma priors are used on all variance components.

```
model
{
  for (i in 1:num.obs)
    {
      X[i,1] <- age[i]
      X[i,2] <- vitAdefic[i]
      X[i,3] <- sex[i]
      X[i,4] <- height[i]
      X[i,5] <- stunted[i]
      X[i,6] <- visit2[i]
      X[i,7] <- visit3[i]
      X[i,8] <- visit4[i]
      X[i,9] <- visit5[i]
      logit(mu[i])
         <- gamma[subject[i]]
         + inprod(beta[],X[i,])
         + inprod(u.spline[],Z.spline[i,])
      resp[i] ~ dbern(mu[i])
    }
  for (i.subj in 1:num.subj)
    {
      gamma[i.subj] <- beta0 + u.subj[i.subj]
      u.subj[i.subj] ~ dnorm(0,tau.u.subj)
    }
  for (k in 1:num.knots)
    {
      u.spline[k] ~ dnorm(0,tau.u.spline)
    }
  beta0 ~ dnorm(0,tau.beta)
  for (j in 1:num.pred)
    {
      beta[j] ~ dnorm(0,tau.beta)
    }
  tau.u.spline ~ dgamma(A.u.spline,
                        B.u.spline)
  tau.u.subj ~ dgamma(A.u.subj,B.u.subj)
}
```

The following code was used to fit (6) to the data on caregiver stress and respiratory health. This code illustrates the use of folded-Cauchy priors on variance components. As noted in the `WinBUGS` user manual (Spiegelhalter, Thomas and Best, 2000), a single zero Poisson observation with mean $\phi$ contributes a term $\exp(\phi)$ to the likelihood for $\sigma$, which is then combined with a flat prior over the positive real line to produce the folded-Cauchy distribution.

```
model
{
  for (i in 1:num.obs)
    {
      X[i,1] <- age[i]
      X[i,2] <- income1[i]
      X[i,3] <- income2[i]
      X[i,4] <- race[i]
```

```
log(mu[i])
   <- gamma[house[i]]
   + inprod(beta[],X[i,])
   + inprod(u.spline[],Z.spline[i,])
y[i] ~ dpois(mu[i])
}
for (i.house in 1:num.house)
{
  gamma[i.house] <- beta0+u.subj[i.house]
  u.subj[i.house] ~ dnorm(0,tau.u.subj)
}
for (k in 1:num.knots)
{
  u.spline[k] ~ dnorm(0,tau.u.spline)
}
beta0 ~ dnorm(0,tau.beta)
for (j in 1:num.pred)
{
  beta[j] ~ dnorm(0,tau.beta)
}
tau.u.spline <- pow(sigma.u.spline,-2)
zero.u.spline <- 0
sigma.u.spline ~ dunif(0,1000)
phi.u.spline <- log((pow(sigma.u.spline,2)
   + pow(phi.scale.u.spline,2)))
zero.u.spline ~ dpois(phi.u.spline)
tau.u.subj <- pow(sigma.u.subj,-2)
zero.u.subj <- 0
sigma.u.subj ~ dunif(0,1000)
phi.u.subj
   <- log((pow(sigma.u.subj,2)
       + pow(phi.scale.u.subj,2)))
zero.u.subj ~ dpois(phi.u.subj)
}
```

Below is the code that we used to fit the spatial model (7) to the Cape Cod female lung cancer data. Please note that the variance components have inverse gamma priors, and `adj`, `weights`, and `num` are inputs to `car.normal`, the normal conditional autoregressive function in `WinBUGS`.

```
model
{
 for (i in 1:num.regions)
  {
   X[i,1] <- working[i,1]
   X[i,2] <- distance[i,1]
   theta[i] <- beta0+u.spatial[i]
       + inprod(beta[],X[i,])
       + inprod(u.spline[],
                Z.spline[i,])
   log(mu[i]) <- log(E[i])+theta[i]
   O[i] ~ dpois(mu[i])
   SIRhat[i] <- 100*mu[i]/E[i]
  }
 u.spatial[1:num.regions]
    ~car.normal(adj[],weights[],
             num[],tau.u.spatial)
for (k in 1:num.knots)
 {
```

```
   u.spline[k] ~ dnorm(0.0,tau.u.spline)
 }
for (j in 1:num.pred)
 {
   beta[j] ~ dnorm(0.0,tau.beta)
 }
beta0 ~ dnorm(0.0,tau.beta)
tau.u.spatial~dgamma(A.u.spatial,
                    B.u.spatial)
tau.u.spline~dgamma(A.u.spline,B.u.spline)
}
```

## ACKNOWLEDGMENTS

We are grateful for comments from Ciprian Crainiceanu, Jim Hodges, Scott Sisson and the reviewers. This research was partially supported by National Institute of Environmental Health Sciences Grant R01-ES10844-01, National Science Foundation Grant DMS-03-06227 and National Institutes of Health Grant ES-01-2044.